\documentclass[a4paper,reqno, 10pt]{amsart}
\usepackage{amsmath,amssymb,amsthm,enumerate}

\usepackage{paralist}
\usepackage{graphics}
\usepackage{epsfig}
\usepackage{amsfonts}
\usepackage{amssymb}
\usepackage{color}

\usepackage{mathtools}

\usepackage[numbers,sort&compress]{natbib}

\usepackage{mathrsfs}
\usepackage{lineno}

\usepackage[svgnames]{xcolor}

\usepackage{pifont}
\usepackage{mathtools}

\usepackage[misc]{ifsym}

\def\ifl{\iffalse }

\def\bc{\begin{center}} \def\ec{\end{center}}
\def\ba{\begin{array}} \def\ea{\end{array}}
\def\bea{\begin{eqnarray}} \def\eea{\end{eqnarray}}
\def\beaa{\begin{eqnarray*}} \def\eeaa{\end{eqnarray*}}

\newtheorem{thm}{Theorem}[section]

\newtheorem{lem}{Lemma}[section]
\theoremstyle{remark}
\newtheorem{rem}{Remark}[section]
\newtheorem*{rem*}{Remark}

\numberwithin{equation}{section}

\newcommand{\R}{\mathbb{R}}

\newcommand{\pa}{\partial}
\newcommand{\na}{\nabla}
\newcommand{\al}{\alpha}
\newcommand{\be}{\beta}

\newcommand{\Ga}{\Gamma}

\renewcommand{\bar}[1]{\overline{#1}}

\newcommand{\Lg}{\langle}
\newcommand{\Rg}{\rangle}

\newcommand{\OK}{\operatorname{OK}}

\title{Uniform boundedness of highest norm for 2D quasilinear
wave}

\author[X.Y. Cheng]{Xinyu Cheng}
\address{X.Y. Cheng, Department of Mathematics, University of British Columbia, Vancouver, BC V6T 1Z2, Canada}
\email{xycheng@math.ubc.ca}

\author[D. Li]{ Dong Li}
\address{D. Li, Department of Mathematics, the Hong Kong University of Science \& Technology, Clear Water Bay, Kowloon, Hong Kong}
\email{mpdongli@gmail.com}

\author[J. Xu]{Jiao Xu}	
\address{J. Xu, SUSTech International Center for Mathematics, Southern University of Science and Technology,
	Shenzhen, P.R. China}
\email{xuj7@sustech.edu.cn}

\begin{document}

\begin{abstract}
We consider the two-dimensional quasilinear wave equations with quadratic nonlinearities. We introduce a new class of null forms and prove uniform boundedness of
the highest order norm of the solution for all time. This class of null forms include several prototypical
strong null conditions as special cases.  To handle the critical decay near the light cone
we inflate the nonlinearity through a new normal form type transformation which is based on a deep
cancelation between the tangential and normal
derivatives with respect to the light cone.  
 Our proof does not employ the 
Lorentz boost and can have promising applications to systems with multiple speeds. 
\end{abstract}

\maketitle

\section{introduction}
 We consider the following quasilinear wave equation:
\begin{align}\label{eq:we2d}
 \begin{cases}
  \square u=g^{kij} \partial_k u \partial_{ij} u,\, \quad t>2, \quad x\in\R^2,\\
  u|_{t=2}=\varepsilon f_1,\ \  \pa_{t} u|_{t=1}=\varepsilon f_2.
 \end{cases}
\end{align}
Here $\square = \partial_{tt} - \Delta$ is the wave operator.  The functions
$f_1$, $f_2$ are real-valued. On the RHS of \eqref{eq:we2d} we
employ the usual Einstein summation convention with $\partial_0 =\partial_t$ and $\partial_l=\partial_{x_l}$ for $l=1,2$.  For simplicity we assume $g^{kij}$ are constant coefficients,
$g^{kij}=g^{kji}$ for any $i$, $j$,  and satisfy the standard null condition:
\begin{align} \label{null1}
g^{kij} \omega_k \omega_i \omega_j=0, 
\end{align}
where $\omega$ is any null vector, namely,  $ \omega_0=-1, \omega_1=\cos\theta,
\omega_2=\sin \theta$ with $\theta \in [0,2\pi]$. In addition, we assume that
\begin{align} \label{null2}
(g^{1ij} \omega_2 - g^{2ij} \omega_1) \omega_i \omega_j =0, \qquad
\text{for any null vector $\omega$}.
\end{align}
We shall call the condition \eqref{null1}--\eqref{null2} a new type of strong null condition. 
Although the extra condition \eqref{null2} seems a bit odd-looking at first sight, we shall
show later that it serves as a  natural generalization of several prototypical null conditions in the literature. 
Our main result reads as follows. 
\begin{thm}\label{thm:main}
Consider \eqref{eq:we2d} with $g^{kij}$ satisfying the strong null condition
\eqref{null1}--\eqref{null2}. Let $m\ge 5$ and assume $f_1 \in H^{m+1}(\mathbb R^2) $, $f_2\in H^{m} (\mathbb R^2)$ are compactly supported in the disk $\{|x| \le 1\}$.  There exists $\varepsilon_{0}>0$ depending on $g^{kij}$ and $\|f_1\|_{H^{m+1}}+\|f_2\|_{H^{m}}$ such that for all $0\le \varepsilon <\varepsilon_0$, the system \eqref{eq:we2d} has a unique global solution. 
Furthermore, the highest norm of the solution remains uniformly bounded, namely
\begin{align}
\sup_{t\ge 2} \sum_{|\alpha|\le m}  \| (\partial \Gamma^{\alpha} u)(t,\cdot) \|_{L_x^2(\mathbb R^2)}
< \infty.
\end{align}
Here $\Gamma=\{ \partial_t, \partial_{x_1}, \partial_{x_2}, \partial_{\theta}, t\partial_t +r\partial_r\}$ does not include the Lorentz boost (see  \eqref{def_Gamma1} for notation). 
\end{thm}

In the seminal work \cite{Alinh01_1},  Alinhac showed that under the general null condition
\eqref{null1} the system \eqref{eq:we2d} has small data global wellposedness with the highest
norm  polynomially bounded in time. In \cite{Lei16}, by deeply exploiting a type of strong
null condition together with Alinhac's method, Lei established small data global wellposedness for 2D incompressible elastodynamics. In \cite{CLM18}, developing upon Lei's strong null condition
in \cite{Lei16}, Cai, Lei and Masmoudi considered the following quasilinear wave equations:
\begin{align}
\square u = A_l \partial_l ( N_{ij} \partial_i u \partial_j u), \label{null3}
\end{align}
where $A_l$, $N_{ij}$ are constants, and 
\begin{align} \label{null3a}
N_{ij} \omega_i \omega_j=0, \qquad\text{for any null vector $\omega$}.
\end{align}
A special case of \eqref{null3} is the following typical quasilinear wave equation
\begin{align} \label{null4}
\square u = \partial_t ( |\partial_t u|^2 - |\nabla u|^2).
\end{align}
In \cite{CLM18} by using a nonlocal transformation (see Remark 1.3 therein) 
it was shown that the system \eqref{null4} has a uniform bound of the highest-order
energy for all time. More recently by using Alinhac's ghost weight
and the null structure in the Lagrangian formulation, Cai \cite{Cai21} showed uniform boundedness of
the highest-order energy for 2D incompressible elastodynamics . In \cite{Dong21}, by using the hyperbolic
foliation method which goes back to H\"ormander and Klainerman, Dong, LeFloch and Lei
showed that the top-order energy of the system \eqref{eq:we2d} with the null condition
\eqref{null1} is uniformly bounded for all time. The main advantage
of the hyperbolic change of variable is that one can gain better control of the conformal
energy thanks to the extra integrability in the hyperbolic time $s=\sqrt{t^2-r^2}$. One should note,
however, that if one works  with the advanced coordinate $s=t-r$, then there is certain degeneracy 
in the $\partial_s$ direction which renders (even any generalized)
conformal energy out of control.  In this connection an interesting further issue is to explore the monotonicity of the conformal energy (and possible generalizations) with respect to different space-time foliations. 

As was already mentioned earlier, the main purpose of this work is to develop a new strategy (building upon Alinhac's ghost weight method) to prove the uniform boundedness of highest norm for general
2D quasilinear wave equations with null conditions. To understand the role of various null
conditions a natural first step is to classify the standard null forms. To this end, we define
the following:
\begin{align}
&F^A_i= \partial_i ( |\partial_t u |^2 - | \nabla u |^2 ), \quad i=0, 1,2; \label{FA00}\\
&F^B_i = \partial_i u \square u, \qquad i=0,1,2; \\
&F^C_1= \partial_0 u \partial_{12} u - \partial_1 u \partial_{02}u; \quad
F^C_2= \partial_1 u \partial_{02} u - \partial_2 u \partial_{01}u;  \label{FA01}\\
&F^D_1= \partial_0 u \partial_{11} u - \partial_1 u \partial_{01} u; \quad
 F^D_2= \partial_1 u \partial_{22} u - \partial_2 u \partial_{12} u; \quad F^D_3= \partial_2 u \partial_{11} u -\partial_1 u \partial_{12} u.  \label{FA11}
\end{align}
\begin{thm}[Full classification of null conditions] \label{thmNull}
If $g^{kij}$  satisfies the standard null condition \eqref{null1}, then
\begin{align}
g^{kij} \partial_k u \partial_{ij} u = \sum_{l=0}^2 C_{1, l} F^A_{l}
+ \sum_{l=0}^2 C_{2,l} F^B_l+ \sum_{l=1}^2 C_{3,l} F^C_l
+\sum_{l=1}^3 C_{4,l} F^D_l,
\end{align}
where $C_{m,l}$ are constants. If $g^{kij}$ satisfies the strong null condition \eqref{null1}--\eqref{null2}, then 
\begin{align} \label{FB001}
g^{kij} \partial_k u \partial_{ij} u = \sum_{l=0}^2 C_{1, l} F^A_{l}
+ \sum_{l=0}^2 C_{2,l} F^B_l. 
\end{align}
If $N_{ij}$ satisfies the null condition \eqref{null3a}, then 
\begin{align} 
A_l \partial_l ( N_{ij} \partial_i u \partial_j u)=\sum_{l=0}^2 C_{1, l} F^A_{l}. 
\end{align}
\end{thm}
\begin{rem}
In yet other words, the null condition in \cite{CLM18} is simply $\partial ( |\partial_t u|^2-|\nabla u|^2)$, whereas our strong null condition is $\operatorname{const} \cdot \partial ( |\partial_tu |^2-|\nabla u|^2)
+\operatorname{const}\cdot \partial u \square u$. In \cite{PZ16}, Peng and Zha 
considered (see formula (1.4) therein) the situation $g^{kij}=g^{ikj}=g^{jik}$ for any
$k$, $i$, $j$ (besides the standard null condition). However, such a strong condition apparently does not include
the standard nonlinearity $\partial ( |\partial_t u|^2 -|\nabla u|^2)$. 
\end{rem}
\begin{rem}
Define the standard null forms
\begin{align}
&G^C_i= \partial_0 u \partial_1 (\partial_i u) - \partial_1 u \partial_0 (\partial_i u), \quad,i=0,1,2;\\
&G^D_i= \partial_1 u \partial_2 (\partial_i u) - \partial_2 u \partial_1(\partial_i u),
\quad i=0,1,2.
\end{align}
It is easy to check that $G^C_2 =F^C_1$, $G^D_0=F^C_2$, $G^C_1=F^D_1$, $G^D_2=F^D_2$,
$G^D_1=-F^D_3$. On the other hand
\begin{align}
G^C_0&= \partial_0u \partial_{01} u -\partial_1 u \Delta u -\partial_1 u \square u  = \frac 12 \partial_1 ( |\partial_t u|^2 -|\partial_1 u|^2) -\partial_1 u \partial_{22} u
-\partial_1 u \square u \notag \\
&= \frac 12 \partial_1 ( |\partial_t u|^2 -|\nabla u|^2) 
+\partial_2 u \partial_{12}u- \partial_1 u \partial_{22} u - \partial_1 u\square u.
\end{align}
Thus $G^C_0$ is a linear combination of $F^A_1$, $F^B_1$ and $F^D_2$. 
\end{rem}

We now explain the key steps of the proof of Theorem \ref{thm:main} (see section 2
for the relevant notation). To elucidate the idea, we fix any multi-index $\alpha$ with
$|\alpha|=m$, and denote $v= \Gamma^{\alpha} u$.  By Lemma \ref{Lem2.3}, we have
\begin{align}
\square v = g^{kij} \partial_k v \partial_{ij} u + \cdots,
\end{align}
where ``$\cdots$" denotes harmless terms which do not contribute to the main term.

Step 1. Weighted energy estimate. We choose $p (r,t)=q(r-t)$ with $q^{\prime}(s)$ nearly
scales as $\langle s \rangle^{-1}$ to derive
\begin{align}
\frac 12 \frac d {dt}
( \| e^{\frac p2} \partial v \|_2^2)
+ \frac 12 \int e^p q^{\prime} \cdot  |Tv|^2 dx =
\int g^{kij} \omega_i \omega_j T_k v  \partial_{t} v \partial_{tt}u e^p dx+\cdots.
\end{align}
The usual strategy is to use Cauchy-Schwartz to derive
\begin{align}
\Bigl|\int g^{kij} \omega_i \omega_j T_k v  \partial_{t} v \partial_{tt}u e^p dx \Bigr|
\le \frac 1 {10} \int |Tv|^2 q^{\prime} e^pdx
+ \operatorname{const}\cdot \int |\partial v|^2 \frac {|\partial_{tt} u|^2} {q^{\prime}} e^p
dx
\end{align}
which yields polynomial growth in time. To resolve this we shall proceed differently.

Step 2. Refined decomposition. By using $T_1=\omega_1\partial_+-\frac{\omega_2}r
\partial_{\theta}$
and $T_2=\omega_2\partial_++\frac{\omega_1} r \partial_{\theta}$,   we have
\begin{align*}
 g^{kij} \omega_{i}\omega_{j}T_{k}v
 &= g^{1ij} \omega_{i}\omega_{j}(\omega_{1}\pa_{+}v-\frac{\omega_{2}}{r}\pa_{\theta}v)
  +g^{2ij} \omega_{i}\omega_{j}(\omega_{2}\pa_{+}v+\frac{\omega_{1}}{r}\pa_{\theta}v)\\
  &= h(\theta) \partial_+v +\omega_i\omega_j (g^{2ij}\omega_1
  -g^{1ij} \omega_2) \frac 1 r \partial_{\theta} v = h(\theta)\partial_+ v,
  \end{align*}
  where  $h(\theta)=g^{1ij} \omega_{1}\omega_{i}\omega_{j}+g^{2ij} \omega_{2}\omega_{i}\omega_{j}$ and we used \eqref{null2} in the last step.

Step 3. Localization, further decomposition and normal form transformation.  We use a bump function
$\phi$ which is localized to $r\sim t$ such that the main piece becomes
\begin{align}
\int h(\theta) \partial_+ v \partial_t v \partial_{tt}u e^p \phi.
\end{align}
The contribution of the regimes $r\le \frac t2$ and $r>2t$ can be shown to be negligible. 
We further use the decomposition $\partial_t = \frac {\partial_++\partial_-}2$ to transform
the main piece as (below we drop the harmless factor $1/2$)
\begin{align} \label{1.19A}
\int h(\theta) \partial_+ v \partial_- v \partial_{tt}u e^p \phi +\text{Negligible}. 
\end{align}
At this point, the crucial observation is to use the fundamental identity 
\begin{align*}
\partial_+ \partial_- =
\square +\frac 1r {\partial_r} +\frac 1 {r^2} {\partial_{\theta\theta}}
\end{align*}
to transform \eqref{1.19A} into an expression which contains an ``inflated" nonlinearity. It is this
novel normal form type transformation which makes the problem subcritical. 

It should be pointed out that we do not employ the usual Lorentz boost vector field in the whole
proof. 
Additionally we developed several new decay estimates for the regime $r\le t/2$ which
was previously un-available due to the lack of Lorentz boost (cf. Lemma \ref{lem2.6}).
Thus this new strategy could have promising applications in systems with multiple speeds. 
\begin{rem}
It is worthwhile pointing out how the symmetry condition $g^{kij}=g^{ikj}=g^{jik}$
for all $k$, $i$, $j$ was needed in \cite{PZ16}. When bounding the quasilinear piece
$\alpha_2=\alpha$, we have
\begin{align}
 \int g^{kij} \partial_k u \partial_{ij} v \partial_t v e^p dx  &=\int g^{kij} \partial_j ( \partial_k u \partial_i v \partial_t v e^p) 
-\int g^{kij} \partial_{kj} u \partial_i v \partial_t v e^p \notag \\
& \qquad -\int g^{kij} \partial_k u \partial_i v \partial_{tj} v e^p
-\int g^{kij} \partial_k u \partial_i v\partial_t v \partial_j (e^p).
\end{align}
By using the symmetry $g^{kij}=g^{kji}$ (this is harmless), we have
\begin{align}
- \int g^{kij} \partial_k u \partial_i v \partial_{tj} v e^p
= -\frac 12 \int g^{kij} \partial_t (\partial_ku \partial_i v \partial_j v e^p)
+\frac 12 \int g^{kij} \partial_{tk} u \partial_i v \partial_j v e^p
+ \frac 12 \int g^{kij} \partial_i v \partial_j v \partial_t (e^p).
\end{align}
Thus 
\begin{align}
 \int g^{kij} \partial_k u \partial_{ij} v \partial_t v e^p dx  &=
 -\int g^{kij} \partial_{kj} u \partial_i v \partial_t v e^p 
 + \frac 12 \int g^{kij} \partial_{tk} u \partial_i v \partial_j v e^p
 + \frac 12 \int g^{kij} \partial_k u 
 \partial_i v \partial_j v \partial_t (e^p)+ \cdots.  \label{intro_rem6a}
 \end{align}
Note that the second term on the RHS of \eqref{intro_rem6a} is not a problem
due to the good decay of $\partial_{tk} u$. On the other hand, in \cite{PZ16} the decay
of $\partial^2 u$ in the regime $r\le t/2$ was not sufficient to treat the first
term on the RHS of \eqref{intro_rem6a}. For this reason (see (3.11) in \cite{PZ16}),
Peng and Zha used the other piece corresponding to $\alpha_1=\alpha$
and the symmetry $g^{kij}=g^{ikj}$ to kill the above term, namely:
\begin{align}
\int g^{kij} \partial_k v \partial_{ij} u \partial_t v e^p
= \int g^{kij} \partial_i v \partial_{kj} u \partial_t v e^p.
\end{align}
One of the main novelty of this work is that we obtained $t^{-\frac 32} \log t$ decay
in the regime $r\le t/2$ (see Lemma \ref{lem2.6}) which can have useful applications
in many other problems.
\end{rem}

The rest of this paper is organized as follows. In Section 2 we collect some preliminaries and useful
lemmas. In Section 3 we give the proof of Theorem \ref{thmNull}. Section 4 and 5 is devoted
to the proof of Theorem \ref{thm:main}.

\subsection*{Acknowledgement.}
D. Li is supported in part by Hong Kong RGC grant GRF 16307317 and 16309518.
We would like to thank Zha Dongbing for some helpful comments.

\section{Preliminaries}

\subsection*{Notation}
We shall us the Japanese bracket notation: $ \langle x \rangle = \sqrt{1+|x|^2}$, for $x \in \mathbb R^d$.  We denote  $\partial_0 = \partial_t$, 
$\partial_i = \partial_{x_i}$, $i=1,2$ and (below $\partial_{\theta}$ and
$\partial_r$ correspond to the usual polar coordinates) 
\begin{align}
& \partial = (\partial_i)_{i=0}^2, \; \partial_{\theta}=x_1 \partial_2 -x_2\partial_1,
\; L_0 = t \partial_t + r \partial_r, \\
& \Gamma= (\Gamma_i )_{i=1}^5, \quad\text{where } \Gamma_1 =\partial_t, \Gamma_2=\partial_1,
\Gamma_3= \partial_2, \Gamma_4= \partial_{\theta}, \Gamma_5=L_0; 
\label{def_Gamma}\\
& \Gamma^{\alpha} =\Gamma_1^{\alpha_1} \Gamma_2^{\alpha_2}
\Gamma_3^{\alpha_3} \Gamma_4^{\alpha_4}\Gamma_5^{\alpha_5}, \qquad \text{$\alpha=(\alpha_1,\cdots, \alpha_5)$ is a multi-index}, \label{def_Gamma1}\\
& \partial_+ =\partial_t + \partial_r, \qquad \partial_- =\partial_t - \partial_r; \\
& T_i = \omega_i \partial_t + \partial_i, \; \omega_0=-1, \; \omega_i=x_i/r, \, i=1,2.
\label{2.5a}
\end{align}
Note that in \eqref{def_Gamma} we do not include the Lorentz boosts. Note that $T_0=0$. 
For simplicity of notation,
we define for any integer $k\ge 1$,  $\Gamma^k = (\Gamma^{\alpha})_{|\alpha|=k}$,
$\Gamma^{\le k} =(\Gamma^{\alpha})_{|\alpha|\le k}$. 
In particular
\begin{align}
|\Gamma^{\le k} u | = (\sum_{|\alpha|\le k} |\Gamma^{\alpha} u |^2)^{\frac 12}.
\end{align}
Informally speaking, it is useful to think of  $\Gamma^{\le k} $ as any one of the vector
fields $ \Gamma^{\alpha}$ with $|\alpha| \le k$. 

For any two quantities $A$, $B\ge 0$, we write  $A\lesssim B$ if $A\le CB$ for some unimportant constant $C>0$. 
%We write $A\lesssim_{Z_1,\cdots,Z_k} B$ if $A\le CB$ where $C>0$ depends on
%the parameters ($Z_1,\cdots, Z_k$).
We write $A\sim B$ if $A\lesssim B$ and $B\lesssim A$. We write $A\ll B$ if
$A\le c B$ and $c>0$ is a sufficiently small constant. The needed smallness is clear from the context.

 \begin{lem}[Sobolev decay]\label{lem:S}
For $v\in C_c^{\infty}(\R^2)$, we have
\begin{equation*}
  | v(x)|\lesssim
  \begin{cases}
    \|v\|_{2}+ \|\Delta v\|_{2}, &\quad |x|\le 1,\\[1.5mm]
    \langle x\rangle^{-\frac12} \|\pa_{\theta}^{\le1}\pa_{r}^{\le 1} v\|_{2}, &\quad |x|>1.
  \end{cases}
\end{equation*}
\end{lem}
\begin{proof}
We focus on the regime $|x|>1$. For a one-variable function $h \in C_c^{\infty}([0,\infty))$, we have
\begin{equation}\label{Sb-1}
    \rho|h(\rho)|^2 \le \int_0^{\infty} |h(r)|^2 r dr + \int_0^{\infty} |\partial_r h|^2 r dr, \quad
    \forall\, \rho>0.
\end{equation}
It follows that (below we slightly abuse the notation and denote $v(\rho,\theta)=v(x)$ for
$x=(\rho\cos\theta, \rho \sin\theta)$)
\begin{align*}
  \rho\|\pa_{\theta}v(\rho,\theta)\|_{L^2_\theta}^2
  =\rho\int_{0}^{2\pi}|\pa_{\theta}v(\rho,\theta)|^2 d\theta
  \lesssim\ \|\pa_{\theta}v \|_{L^2(\mathbb R^2)}^2+\|\pa_{r}\pa_{\theta}v \|_{L^2(\mathbb R^2)}^2.
\end{align*}
Denote $\bar v(\rho)$ as the average of $v(\rho,\theta)$ over $\theta$. By \eqref{Sb-1}, we have
\begin{align*}
  \rho|\bar{v}(\rho)|^2
  ={\rho}\left( \frac 1 {2\pi}\int_{0}^{2\pi}v(\rho,\theta) d\theta\right)^2
\lesssim \rho \int_0^{2\pi} v(\rho,\theta)^2 d\theta
  \lesssim \, \|v \|_{L^2(\mathbb R^2) }^2+\|\pa_{r}v \|_{L^2(\mathbb R^2)}^2.
\end{align*}

Note that $|v(\rho,\theta)-\bar{v} (\rho)|^2\lesssim|\pa_{\theta}v|_{L^2_{\theta}}^2$ by the Poincar\'e inequality. Thus
\begin{equation*}
  |x||v(x)|^2
  \lesssim\, \|v \|_{2}^2+ \|\pa_{r}v\|_{2}^2
  +\|\pa_{\theta}v\|_{2}^2+\|\pa_{r}\pa_{\theta}v\|_{2}^2.
\end{equation*}
\end{proof}

\begin{lem}[Refined Hardy's inequality]\label{lem:Hardy}
For any real-valued $h\in C_c^{\infty}([0, M+1))$ with $M>0$, we have
\begin{align} \label{2.6.0a}
\int_0^{M+1} \frac {h(\rho)^2} {(2+M-\rho)^2} \rho d\rho 
\le  \, 4\int_0^{\infty} (h^{\prime}(\rho) )^2 \rho d \rho.
\end{align}

For $u\in C^{\infty}([0,T]\times \R^2)$ with support in $\{(t,x): |x|\leq 1+t\}$, we have
\begin{align*}
 \|\langle |x|-t\rangle^{-1} u \|_{L_x^2(\mathbb R^2) }&\lesssim \|\pa_r u\|_{L_x^2(\mathbb R^2)},
  \qquad     \langle |x|-t\rangle^{-1} |u(t,x)|\lesssim 
  \langle x\rangle^{-\frac 12}  \| \partial \Gamma^{\le 1} u \|_{L_x^2(\mathbb R^2)}. 
\end{align*}

\end{lem}

\begin{proof}
The inequality \eqref{2.6.0a} follows from integrating by parts:
\begin{align}
\text{LHS of \eqref{2.6.0a}}= - \int_0^{M+1}
\frac{h^2} {2+M-\rho} d\rho +\int_0^{M+1} \frac {2hh^{\prime}}{2+M-\rho} \rho d\rho.
\end{align}
The second inequality follows from \eqref{2.6.0a} and the 
fact that $\langle |x|-t\rangle^{-2} \sim (2+t-|x|)^{-2}$ for $|x|\le 1+t$. 
For the third inequality, consider first the case $|x|>1$.  By Lemma \ref{lem:S}, we have
\begin{align}
\langle |x|-t \rangle^{-1} |u(t,x)|
& \lesssim \langle x \rangle^{-\frac 12}  \| \partial_r^{\le 1}
\partial_{\theta}^{\le 1} ( \langle r-t\rangle^{-1} u ) \|_2 \sim
\langle x \rangle^{-\frac 12} | \partial_r^{\le 1} ( \langle r -t \rangle^{-1} \partial_{\theta}^{\le 1} 
u) \|_2 \notag \\
& \lesssim \langle x \rangle^{-\frac 12} \| \partial_r \partial_{\theta}^{\le 1} u\|_2 
\lesssim \langle x \rangle^{-\frac 12}
 \| \partial \Gamma^{\le 1} u \|_2.
\end{align}
On the other hand, for $|x|\le 1$, we have
\begin{align*}
\langle |x| - t \rangle^{-1} |u(t,x)| 
 & \lesssim \langle t \rangle^{-1}  ( \| u \|_{L^2_x(|x| \le 1)} + \| \partial^2 u \|_{L_x^2(|x|\le 1)} )  \notag \\
 & \lesssim \| \langle |x|- t\rangle^{-1} u \|_{L_x^2(\mathbb R^2)}+ \| \Delta u \|_{L_x^2(\mathbb R^2)} \lesssim \| \nabla u \|_2
 +\|\Delta u\|_2.
 \end{align*}
Note that we actually proved the inequality: 
\begin{align*}
 \langle x\rangle^{\frac12}\langle |x|-t\rangle^{-1} |u(t,x)|\lesssim \begin{cases}
 \|\pa_r\pa_{\theta}^{\le 1} u \|_{2},\quad \ &|x|\ge 1;\\
    \|\nabla u \|_{2}+\|\Delta u\|_{2},\quad\ &|x|<1.
     \end{cases}
\end{align*}
\end{proof}

\begin{lem} \label{Lem2.3}
If $g^{kij}$ satisfies the null condition, then for $t> 0$  we have
\begin{align} \label{2.10A}
g^{kij} \partial_k f
\partial_{ij} h   = g^{kij}
(T_k  f \partial_{ij} h
-\omega_k \partial_t f T_i \partial_j h
+ \omega_k \omega_i \partial_t f T_j
\partial_t h),
\end{align}
where $T=(T_1,T_2)$ is defined in \eqref{2.5a}.  It follows that
\begin{align}
|g^{kij} \partial_k f \partial_{ij} h| 
& \lesssim | T f | |\partial^2 h| + |\partial f | | T \partial h|  \label{a2.12a}\\
& \lesssim 
\frac 1 {\langle r +t \rangle} 
(|\Gamma f| |\partial^2 h| + |\partial f | |\Gamma \partial h|
+ |\partial f | \cdot |\partial^2 h| \cdot |r -t | ). \label{a2.12b}
\end{align}
Suppose $g^{kij}$ satisfies the null condition and 
$
\square u = g^{kij} \partial_k u \partial_{ij} u.
$
Then for any multi-index $\alpha$, we have
\begin{align} \label{a2.12aa}
\square \Gamma^{\alpha} u  = \sum_{\alpha_1+ \alpha_2 \le \alpha}
g^{kij}_{\alpha_1,\alpha_2} \partial_k \Gamma^{\alpha_1} u 
\partial_{ij} \Gamma^{\alpha_2} u,
\end{align}
where for each ($\alpha_1$, $\alpha_2$), $g^{kij}_{\alpha_1,\alpha_2}$ also satisfies the
null condition.  In addition, we have $g^{kij}_{\alpha,0} =g^{kij}_{0,\alpha}=g^{kij}$.
\end{lem}
\begin{proof}
The identity \eqref{2.10A} follows by applying repeatedly the identity $\partial_l
=T_l - \omega_l \partial_t$ and using the null condition at the last step. 
The inequality \eqref{a2.12b} is obvious if $r\le \frac t2$ or $r\ge 2t$, or $r \sim t \lesssim 1$ since
$\langle r +t \rangle \sim \langle r-t \rangle $ in these regimes. On the other hand, if
$r\sim t\gtrsim 1$,  then one can use the identities
\begin{align}
&T_1 = \omega_1 \partial_+ -\frac {\omega_2} r \partial_\theta, \quad
T_2 = \omega_2 \partial_+ +\frac {\omega_1} r \partial_\theta ;\quad
 \partial_+ = \frac 1 {t+r} ( 2 L_0 - (t-r) \partial_-).
\end{align}
The identity \eqref{a2.12aa} follows from H\"ormander \cite{Hor97}.
\end{proof}

\begin{lem} \label{lem2.3a}
Suppose $\tilde u= \tilde u(t,x)$ has continuous second order derivatives. Then
\begin{align} 
&| \langle r -t \rangle \partial_{tt} \tilde u (t,x) | 
+| \langle r -t \rangle \partial_t \nabla  \tilde u (t,x) | 
+| \langle r -t \rangle \Delta \tilde u (t,x) |  \notag \\
\lesssim &
| (\partial \Gamma^{\le 1} \tilde u)(t,x)| + (r+t) | (\square \tilde u)(t,x) |, \quad r=|x|, \, 
t\ge 0;  \label{2.9a0}
\end{align}
and 
\begin{align}
&| \langle r -t \rangle  \partial^2 \tilde u (t,x) | 
\lesssim 
| (\partial \Gamma^{\le 1} \tilde u)(t,x)| + (r+t) | (\square \tilde u)(t,x) |, \quad 
\forall\, r\ge t/10,
 \, t\ge 1. \label{2.9a1}
\end{align}
Suppose $T_0\ge 1$ and $u \in C^{\infty}([1,T_0]\times \mathbb R^2$ solves  \eqref{eq:we2d} with support in $|x|\le t+1$, $1\le t\le T_0$.
For any integer $l_0\ge 2$, there exists $\epsilon_1>0$ depending only on
$l_0$, such that if at some $1\le t\le T_0$,
\begin{align} \label{2.9a3}
\| (\partial \Gamma^{\le \lceil{\frac {l_0}2}\rceil +2} u)(t,\cdot) \|_{L_x^2(\mathbb R^2)} \le \epsilon_1,
\qquad  \text{( here $\lceil{z}\rceil = \min \{n\in \mathbb N:\, n\ge z\}$ )}
\end{align}
then for the same $t$, we have the $L^2$ estimate:
\begin{align} \label{2.9a4}
\| (\langle r -t \rangle \partial^2 \Gamma^{\le l_0} u) (t,\cdot) \|_{L_x^2(\mathbb R^2)}\lesssim \| (\partial \Gamma^{\le l_0+1} u )(t,\cdot) \|_{L_x^2(\mathbb R^2)}.
\end{align}
For any integer $l_1\ge 0$, there exists $\epsilon_2>0$ depending only on
$l_1$, such that if at some $1\le t\le T_0$,
\begin{align} \label{2.99a1}
\| (\partial \Gamma^{\le l_1 +2} u)(t,\cdot) \|_{L_x^2(\mathbb R^2)} \le \epsilon_2,
\end{align}
then for the same $t$, we have the point-wise estimate:
\begin{align} \label{2.99a2}
 |(\langle r -t  \rangle \partial^2 \Gamma^{\le l_1} u )(t,x) |
 \lesssim |  (\partial \Gamma^{\le l_1+1} u )(t,x)|, 
 \qquad\forall\, r\ge t/10.
 \end{align}
 \end{lem}
\begin{proof}
In the 3D case, the estimate \eqref{2.9a0} is an elementary but deep observation of Sideris (cf. \cite{SK96}). Note that for 2D by using $\partial_{\theta} L_0 \tilde u = (t-r) \partial_{\theta}
\partial_t \tilde u + r (\partial_t+\partial_r) \partial_{\theta} \tilde u$, we have
$
|t-r|\cdot |\frac 1 r \partial_{\theta} \partial_t \tilde u | \lesssim |\partial \Gamma^{\le 1} \tilde u|
$
which (together with the estimate of $|\partial_t \partial_r \tilde u|$) 
settles the estimate for $|\nabla \partial_t u|$. 
The estimate \eqref{2.9a1} can be derived along similar lines since $r\ge t/10$ (note that $\partial^2$ includes
$\partial_i \partial_j$!). 
%In the appendix we give an intuitive proof for the sake of completeness. 
 For \eqref{2.9a4}, by using a simple integration-by-parts argument, one has (below 
 $k_0\ge 0$ is a  running parameter)
 \begin{align} \label{2.21A}
\sum_{i,j=1}^2 \| \langle r -t\rangle \partial_i \partial_j \Gamma^{\le k_0} u \|_2
  \lesssim \| \partial \Gamma^{\le k_0} u \|_2 + \| \langle r -t\rangle \Delta 
  \Gamma^{\le k_0} u\|_2.
  \end{align}
 By using \eqref{2.9a0} and
\eqref{a2.12b} we have
\begin{align}
 & |(\langle r -t \rangle \partial_{tt} \Gamma^{\le k_0} u)(t,x)|  +
 |(\langle r -t \rangle \partial_t \nabla  \Gamma^{\le k_0} u)(t,x)|  +|(\langle r -t \rangle \Delta \Gamma^{\le k_0} u)(t,x)|   \notag \\
\lesssim  &
| (\partial \Gamma^{\le k_0+1}  u)(t,x)| 
+ \sum_{m+l\le k_0}
(|\Gamma^{\le m+1} u|
|\partial^2 \Gamma^{\le l} u |
+|\partial \Gamma^{\le m} u| |\Gamma^{\le l+1} \partial u|+
|\partial \Gamma^{\le m} u | |\partial^2 \Gamma^{\le l} u| |r-t| ). \label{2.23P}
\end{align}
By \eqref{2.21A}, we obtain
\begin{align}
&\| \langle r-t\rangle \partial^2 \Gamma^{\le k_0} u \|_2 
\lesssim  
\| \partial \Gamma^{\le k_0+1}  u\|_2  \notag \\
&\qquad+ \sum_{m+l\le k_0}
(\| |\Gamma^{\le m+1} u|
|\partial^2 \Gamma^{\le l} u | \|_2
+\| |\partial \Gamma^{\le m} u| |\Gamma^{\le l+1} \partial u| \|_2+
\| |\partial \Gamma^{\le m} u | |\partial^2 \Gamma^{\le l} u| |r-t| \|_2).
\end{align}
If $m\le l+1$, then we use  the estimates 
(note that $m+2 \le \lfloor
\frac{k_0+1}2 \rfloor +2\le \lceil \frac {k_0}2\rceil+2$)
%By Lemma \ref{lem:Hardy}, we have 
\begin{align}
&\langle r -t\rangle^{-1}| (\Gamma^{\le m+1} u)(t,x)| \lesssim
\| \partial \Gamma^{\le m+2} u \|_2, 
\quad \| \partial \Gamma^{\le m} u\|_{\infty}
\lesssim \| \partial \Gamma^{\le m+2} u\|_2.
\end{align}
If $m\ge l+2 $, then $l\le \frac{k_0-2}2$ and we use the estimates 
(see \eqref{2.26B} for the second estimate)
\begin{align} \label{2.25A}
\| \frac{|\Gamma^{\le m+1} u|}{\langle r-t\rangle} \|_2 \lesssim
\| \partial \Gamma^{\le m+1} u\|_2, \quad
| \langle r -t \rangle \partial^2 \Gamma^{\le l} u(t,x)|
\lesssim \|\langle r -t\rangle \partial^2 \Gamma^{\le l+2} u\|_2.
\end{align}
Thus if $\| \partial\Gamma^{\le k_0+1} u \|_2 \ll 1$, we obtain 
\begin{align}
\| \langle r-t \rangle \partial^2 \Gamma^{\le k_0} u(t,\cdot) \|_2
\lesssim \| \partial \Gamma^{\le k_0+1} u\|_2 \ll 1.
\end{align}
To prove \eqref{2.9a4} under the assumption \eqref{2.9a3} we first take $k_0=
\lceil{\frac {l_0}2}\rceil+1 $ and show that 
\begin{align}
\|\langle r-t\rangle \partial^2 \Gamma^{\le \lceil{\frac {l_0}2}\rceil +1} u)(t,\cdot) \|_2 
\lesssim\| (\partial \Gamma^{\le \lceil{\frac {l_0}2}\rceil +2} u)(t,\cdot) \|_2 \ll 1.
\end{align}
We then use this smallness in \eqref{2.25A} and obtain the desired result for
$k_0= l_0$ (Note that $\lceil \frac{l_0-2} 2 \rceil +2 \le \lceil \frac{l_0} 2 \rceil+1$). 
The estimate of \eqref{2.99a2} follows from  \eqref{2.9a1}.
\end{proof}

\begin{lem} \label{lem2.5A}
For any $f \in C_c^{\infty}(\mathbb R^2)$, we have
\begin{align} 
&\langle |x_0| -t \rangle^{\frac 12} |f(x_0)| \lesssim \|f\|_2+\| \langle |x|-t\rangle \nabla f
\|_2 + \| \langle |x|-t\rangle \partial_1 \partial_2 f\|_2, \quad \forall\, x_0\in \mathbb R^2,
\, t\ge 0;  \label{2.26A}\\
& \| \langle |x| -t \rangle \partial f \|_{\infty} 
\lesssim \| \langle |x|-t\rangle \partial f \|_2 + \| \langle |x| -t\rangle \partial^2 f \|_2
+ \| \langle |x| -t \rangle \partial^3 f \|_2, \, \quad\forall\, t\ge 0.
\label{2.26B}
\end{align}
It follows that
\begin{align} 
&\| f \|_{L_x^{\infty}(\mathbb R^2)} \lesssim
\langle t \rangle^{-\frac 12}
(\|f\|_2 + \| \langle |x| -t\rangle \nabla \tilde \Gamma^{\le 1 } f \|_2), \qquad\forall\, t\ge 0, 
\label{2.27A}
%&\| \partial f \|_{L_x^{\infty}(|x|<\frac 23 t)}
%\lesssim \langle t\rangle^{-1} \| \langle |x|-t\rangle 
%\partial \Gamma^{\le 2} f \|_2, \qquad \forall\, t\ge 0.
\end{align}
where $\tilde \Gamma =(\partial_1, \partial_2, \partial_{\theta})$. 
\end{lem}
\begin{proof}
The case $||x_0|-t|\le 2$ follows from the inequality
$|f(x_0)|^2 \le \int|\partial_1\partial_2(f(x)^2) | dx_1 dx_2$. For $||x_0|-t|>2$,
we note that $\langle |x_0|-t \rangle \lesssim \frac{ \langle |x_0|^2-t^2 \rangle}{\langle x_0
\rangle + t}$.   
We need to work with the latter weight since it will be smooth near the spatial
origin. Note that $\frac{\langle |x|^2 -t^2 \rangle}{\langle x \rangle +t}
\lesssim \langle |x|-t \rangle $ for all $x\in \mathbb R^2$, $t\ge 0$. This will be used
in the computation below.

Since $f \in C_c^{\infty}$, by using the Fundamental Theorem of Calculus we have 
\begin{align}
 \langle |x_0|-t \rangle | f(x_0)|^2 &\lesssim  \frac{\langle |x_0|^2 -t^2  \rangle} {\langle x_0
\rangle +t} |f(x_0)|^2 \lesssim \int_{\mathbb R^2}
\left|\partial_1\partial_2 \Bigl ( \frac{\langle |x|^2 -t^2 \rangle} {\langle x \rangle + t} f(x)^2
\Bigr) \right | dx_1 dx_2 \notag \\
& \lesssim \| f\|_2^2 + \| \nabla f \|_2 \| f\|_2 + \|
 \langle |x|-t\rangle \nabla f \|_2 \| \nabla f\|_2  + \| \langle |x| -t\rangle \partial_1\partial_2 f \|_2
 \| f\|_2\notag \\
 & \lesssim \|f\|_2^2 + \| \langle |x|-t\rangle \nabla f \|_2^2+ \|
 \langle |x|-t\rangle \partial_1 \partial_2 f \|_2^2.
\end{align}
Thus \eqref{2.26A} follows. The proof of \eqref{2.26B} is similar.  Let $W(x) = \frac{\langle |x|^2 -t^2 \rangle}{\langle x \rangle +t}$ and observe that
\begin{align}
\sum_{1\le i\le 2}\|\partial_i  W\|_{\infty}+ \sum_{1\le i,j\le 2}\|\partial_i \partial_j W\|_{\infty} \lesssim 1.
\end{align}
One can then work with the expression $W(x_0)^2 |\partial f(x_0)|^2$ for the case $|x_0-t|>2$
and derive the result.
 For \eqref{2.27A} we may assume $t\ge 2$. The
case $|x_0|\le t/2$ follows from \eqref{2.26A}. The case $|x_0|>t/2$ follows from
Lemma \ref{lem:S}. 
\end{proof}

\begin{lem}[Decay estimates] \label{lem2.6}
Suppose $T_0\ge 2$ and $u \in C^{\infty}([2,T_0]\times \mathbb R^2)$ solves  \eqref{eq:we2d} with support in $|x|\le t+1$, $2\le t\le T_0$. Suppose $m\ge 4$ and 
\begin{align}
E_m(u(t,\cdot) ) = \| (\partial \Gamma^{\le m} u )(t,\cdot) \|_2^2 \le \tilde \epsilon,
\end{align}
where $\tilde \epsilon>0$ is sufficiently small. Then we have the following decay estimates:
\begin{align}
& t^{\frac 12} \| \partial \Gamma^{\le m-2} u \|_{L_x^{\infty}} 
+t^{\frac 12}\| \langle |x|-t\rangle \partial^2 \Gamma^{\le m-3} u \|_{L_x^{\infty}(|x|>\frac{t}{10})} 
+\| \langle |x|-t\rangle \partial^2 \Gamma^{\le m-3} u \|_{L_x^{\infty}}
\lesssim \;1; \label{2.30A}\\
& \| \Delta \Gamma^{\le m-3} u \|_{L_x^{\infty}(|x|<\frac {2t}3)}
+\| \partial_{tt}\Gamma^{\le m-3} u \|_{L_x^{\infty}(|x|<\frac {2t}3)}
+\| \partial_t \nabla  \Gamma^{\le m-3} u \|_{L_x^{\infty}(|x|<\frac {2t}3)}
\lesssim t^{-\frac 32}; \qquad \label{2.30AA}\\
& \sum_{i,j=1}^2 \| \partial_i \partial_j \Gamma^{\le m-3} u
\|_{L_x^{\infty}(|x| <\frac t2)} \lesssim t^{-\frac 32}
\log t, \quad
\| \langle |x|-t\rangle \partial^2 \Gamma^{\le m-3} u \|_{L_x^{\infty}}
\lesssim t^{-\frac 12} \log t;
\qquad \label{2.30AB}\\
&\| \frac{T \Gamma^{\le m-2} u} {\langle |x|-t \rangle} \|_{L_x^{\infty}} 
+\| T \partial \Gamma^{\le m-3} u \|_{L_x^{\infty}(|x|>\frac t2)}
\lesssim t^{-\frac 32}, \quad \| T \partial \Gamma^{\le m-3} u \|_{L_x^{\infty}} \lesssim 
\langle t \rangle^{-\frac 32} \log t;  \label{2.30B} \\
& \| \frac{ T \Gamma^{\le m-1} u } { \langle |x| -t \rangle }
\|_{L_x^2}  +\| T \partial \Gamma^{\le m-1} u \|_{L_x^2} \lesssim 
t^{-1}; \label{2.30C} \\
& \| \langle |x|-t\rangle^2 \partial^3 u \|_{L_x^{\infty}(|x|\ge \frac t2)}
\lesssim t^{-\frac 12}. \label{2.30D}
\end{align}
\end{lem}
\begin{proof}
We shall take $\tilde \epsilon$ sufficiently small so that Lemma
\ref{lem2.3a} can be applied.  The estimate \eqref{2.30A} follows from Lemma \ref{lem2.5A} and Lemma \ref{lem2.3a}.
The estimate \eqref{2.30AA} follows from \eqref{2.23P}.
Note that the term $ \|(\Gamma^{\le m+1} u)(t) \|_{\infty}$ therein can be bounded using
a dispersive estimate which yields $t^{\frac 12}$ growth.
The estimate \eqref{2.30AB} follows from a harmonic analysis estimate
using \eqref{2.30AA} and separation of supports. 
For \eqref{2.30B}, we note that the case $|x|\le \frac t2$ follows from 
\eqref{2.30A}--\eqref{2.30AB}:
\begin{align}
& \|\frac{T \Gamma^{\le m-2} u} {\langle |x|-t \rangle} \|_{L_x^{\infty}(|x|\le \frac t2)} 
\lesssim  t^{-1} \| \partial \Gamma^{\le m-2} u \|_{L_x^{\infty} (|x|\le \frac t2)}
\lesssim t^{-\frac 32}; \\
& \| T \partial \Gamma^{\le m-3} u \|_{L_x^{\infty}(|x|\le \frac t2)} 
\lesssim t^{-\frac 32} \log t.
\end{align}
On the other hand, for $|x|>\frac t2$ we denote $\tilde u=\Gamma^{\le m-2} u$ and estimate
$\| \frac {T_1 \tilde u} {\langle |x|-t \rangle} \|_{L_x^{\infty}(|x|>\frac t2)}$ (the
estimate for $T_2$ is similar). Recall that
\begin{align}
T_1 \tilde u & = \omega_1 \partial_t \tilde u + \partial_1 \tilde u
= \omega_1 (\partial_t + \partial_r ) \tilde u- \frac{\omega_2} r \partial_{\theta} \tilde u\notag \\
& = \omega_1 \frac 1 {t+r} ( 2 L_0 \tilde u - (t-r) \partial_-\tilde u) -\frac{\omega_2} r \partial_{\theta}
\tilde u.
\end{align}
Clearly for $r=|x|\ge \frac t2$, 
\begin{align}
\left| \frac{ T_1 \tilde u } { \langle r -t \rangle}
\right| &\lesssim \frac1 t \Bigl( \left| \frac{L_0 \tilde u} {\langle r -t\rangle } \right| + |\partial \tilde u|
\Bigr)
+ \left| \frac{\partial_{\theta} \tilde u} {r \langle r-t \rangle } \right| \notag \\
& \lesssim t^{-1} \cdot t^{-\frac 12} \| \partial \Gamma^{\le 2} \tilde u\|_2
+ t^{-\frac 32}+t^{-1} \cdot t^{-\frac 12} \| \partial \Gamma^{\le 1} \partial_{\theta} \tilde u\|_2 \lesssim
t^{-\frac 32},
\end{align}
where in the second last step we used Lemma \ref{lem:Hardy} (for the term $|\partial \tilde u|$
we use \eqref{2.30A}). The estimates for other terms 
in \eqref{2.30AA}--\eqref{2.30C}  are similar. We now sketch how to
prove \eqref{2.30D}.  By using \eqref{2.9a1} (applied to $\tilde u =\partial u$),
we obtain
\begin{align}
|\langle r -t\rangle \partial^3 u|
\lesssim | \partial^2 \Gamma^{\le 1} u |
+ (r+t) | \square \partial u |.
\end{align}
The contribution of the term $|\partial^2 \Gamma^{\le 1} u|$ is clearly OK for us since it can
absorb a factor of $\langle r -t \rangle$. By Lemma \ref{Lem2.3}, we have
\begin{align}
(r+t) |\square \partial u |
&\lesssim
|\Gamma \partial u | |\partial^2 u | + |\partial^2 u |^2 |r-t|
+|\Gamma u| |\partial^3 u | + |\partial u | |\Gamma \partial^2 u|
+ |\partial u | |\partial^3 u | |r-t| \notag \\
& \lesssim |\partial \Gamma^{\le 1} u |
| \partial^2 \Gamma^{\le 1} u|+|\partial^2 u |^2 |r-t|
+|\frac{\Gamma u}{\langle r-t\rangle}| \cdot  \langle r-t\rangle |\partial^3 u |
+ |\partial u | |\partial^3 u | |r-t|.
\end{align}
The desired estimate then clearly follows by using smallness of the pre-factors.
\end{proof}

\section{Proof of Theorem \ref{thmNull} }
Denote (below $w_0=-1$, $w_1=\cos \theta$, $w_2=\sin \theta$)
\begin{align}
\alpha_k &=g^{kij} w_i w_j \\
& = g^{k00} + g^{k01} (-2) \cos \theta + g^{k02} (-2) \sin \theta
+g^{k11} \cos^2 \theta + g^{k22} \sin^2 \theta +g^{k12} 2\cos \theta \sin \theta \\
& =(g^{k00} + \frac {g^{k11}+g^{k22}}2) 
-2g^{k01} \cos \theta -2g^{k02} \sin \theta
+\frac 12 (g^{k11} - g^{k22}) \cos 2\theta +g^{k12} \sin 2\theta.
\end{align}

We compute
\begin{align}
&\alpha_1 w_1  \notag \\
=&(g^{100} + \frac {g^{111}+g^{122}}2) \cos \theta
-2g^{101} \cos^2 \theta -g^{102} \sin 2\theta
+\frac  {g^{111} - g^{122}}2 \cos 2\theta\cos \theta +g^{112} \sin 2\theta \cos \theta \\
 =& - g^{101}+(g^{100} + \frac {g^{111}+g^{122}}2 +\frac  {g^{111} - g^{122}}4) \cos \theta + \frac {g^{112}} 2\sin \theta -{g^{101}} \cos 2\theta -g^{102} \sin 2\theta
 \notag \\
 &\qquad +(\frac  {g^{111} - g^{122}}4) \cos3 \theta+\frac {g^{112}} 2\sin 3\theta.
\end{align}
\begin{align}
& \alpha_2 w_2 \notag \\
=& (g^{200} + \frac {g^{211}+g^{222}}2)  \sin \theta
-g^{201} \sin 2 \theta -2g^{202} \sin^2 \theta
+\frac 12 (g^{211} - g^{222}) \sin \theta \cos 2\theta +g^{212} \sin 2\theta \sin \theta
\notag \\
= & - {g^{202}}  +  \frac {g^{212}} 2 \cos \theta+
 (g^{200} + \frac {g^{211}+g^{222}}2-
 \frac {g^{211} - g^{222}}4)  \sin \theta+g^{202} \cos 2\theta-g^{201} \sin 2 \theta \notag \\
 & \qquad -\frac{g^{212}}2 \cos 3\theta+(\frac {g^{211} - g^{222}}4) \sin 3\theta  .
\end{align}

From the identity $\alpha_0= \alpha_1 w_1 + \alpha_2 w_2$, we obtain
\begin{align}
& g^{101}=-\frac{g^{000}+g^{011}}2, \quad g^{202}=-\frac{g^{000}+g^{022}}2, \\
&  -2 g^{001} = g^{100} + g^{111}  ; \quad -2g^{002} =g^{200} + g^{222}; \\
 &g^{012} =-g^{102}-g^{201}; \quad g^{212}=\frac  {g^{111} - g^{122}}2; \quad  {g^{112}} =-\frac {g^{211} - g^{222}}2.
\end{align}
We now set (below we shall treat $(a_j)_{j=1}^{11}$ as free parameters)
\begin{align}
&g^{000}=a_1, \; g^{001} =a_2, \; g^{002} =a_3, \;
g^{011}=a_4, \; g^{012}=a_5, \; g^{022} =a_6, \\
&g^{100}=a_7, \; g^{222}=a_8, \;g^{102}=a_9, \; g^{122}=a_{10},\;
g^{112}=a_{11}. 
\end{align}
Then we obtain the rest of the coefficients as follows:
\begin{align}
&g^{101}=-\frac{a_1+a_4}2, \;  g^{202}= -\frac{a_1+a_6}2,\;
g^{111}=-2a_2-a_7, \\
&g^{200}=-2a_3-a_8, \, g^{201}=-a_5-a_9,\,
g^{212}=- \frac{2a_2+a_7+a_{10}}2, \\
& g^{211}= a_8-2a_{11}.
\end{align}
We then rearrange the nonlinearity $g^{kij} \partial_k u \partial_{ij} u$ as 11 terms:
\begin{itemize}
\item $a_1$:  $\partial_0 u \partial_{00} u - \partial_1 u \partial_{01} u -
\partial_2 u \partial_{02} u$. 
\item $\frac 12a_2$:  $\partial_0 u \partial_{01} u - \partial_1 u \partial_{11} u
- \partial_2 u \partial_{12} u$. 
\item $\frac 12 a_3$:  $\partial_0 u \partial_{02} u -\partial_2 u \partial_{00} u$.
\item $a_4$: $\partial_0 \partial_{11} u -\partial_1 u \partial_{01} u$.
\item $\frac 12a_5$:  $\partial_0 \partial_{12} u-\partial_2 u \partial_{01} u$.
\item $a_6$: $\partial_0 u \partial_{22} u - \partial_2 u \partial_{02} u$.
\item $a_7$: $\partial_1 u \partial_{00} u - \partial_1 u \partial_{11}u -\partial_2 u \partial_{12} u$.
\item $a_8$: $\partial_2u \partial_{22} u -\partial_2 u \partial_{00} u + \partial_2 u \partial_{11} u$.
\item $\frac 12 a_9$: $\partial_1 u \partial_{02} u -\partial_2 u \partial_{01} u$.
\item $a_{10}$: $\partial_1 u \partial_{22}u-\partial_2 u \partial_{12} u$.
\item $\frac 12a_{11}$: $\partial_1 u \partial_{12} u - \partial_2 u \partial_{11} u$.
\end{itemize}
It is not difficult to check that these 11 terms are in one-to-one correspondence
of \eqref{FA00}--\eqref{FA11} (after suitable linear combinations).

\subsection{Our new null condition}
Now consider the expression 
\begin{align}
& (-\sin \theta) \alpha_1 + (\cos \theta) \alpha_2 \notag \\
= &(-\sin \theta) \cdot \Bigl( 
g^{100}+\frac{g^{111} +g^{122}}2 -2g^{101} \cos \theta
-2 g^{102} \sin \theta+
\frac {g^{111}-g^{122}} 2 \cos 2\theta 
+ g^{112} \sin 2\theta\Bigr) \notag \\
& \; + (\cos \theta) \Bigl( 
g^{200}+\frac{g^{211} +g^{222}}2 -2g^{201} \cos \theta
-2 g^{202} \sin \theta+
\frac {g^{211}-g^{222}} 2 \cos 2\theta 
+ g^{212} \sin 2\theta\Bigr) \notag \\
=& g^{102}-g^{201} + \Bigl(-\frac{g^{112}}2+g^{200}+\frac{g^{211}+g^{222}}2
+\frac{g^{211}-g^{222}}4 \Bigr)\cos \theta \notag \\
& \quad + \Bigl( - (g^{100}+\frac{g^{111}+g^{122}}2)
+\frac{g^{111}-g^{122}}4 +\frac{g^{212}}2 ) \sin \theta \notag \\
& + (-g^{102}-g^{201} ) \cos 2\theta +  (g^{101}-g^{202}  ) \sin 2\theta \notag \\
& +\Bigl( \frac{g^{112}}2+ \frac{g^{211}-g^{222}}4 \Bigr) \cos 3\theta + 
\Bigl(  -\frac{g^{111}-g^{122}}4 +\frac{g^{212}}2\Bigr) \sin 3\theta.
\end{align}
Simplifying a bit using the standard null condition, we obtain
\begin{align}
 & (-\sin \theta) \alpha_1 + (\cos \theta) \alpha_2 \notag \\
 =& g^{102}-g^{201} + \Bigl(-\frac{g^{112}}2+g^{200}+\frac{g^{211}+g^{222}}2
+\frac{g^{211}-g^{222}}4 \Bigr)\cos \theta \notag \\
& \quad + \Bigl( - (g^{100}+\frac{g^{111}+g^{122}}2)
+\frac{g^{111}-g^{122}}4 +\frac{g^{212}}2 ) \sin \theta \notag \\
& + (-g^{102}-g^{201} ) \cos 2\theta +  (g^{101}-g^{202}  ) \sin 2\theta. 
\end{align}
Forcing $  (-\sin \theta) \alpha_1 + (\cos \theta) \alpha_2 \equiv 0$ and using the standard
null condition gives us
\begin{align}
&g^{102}=g^{201} =0, 
\quad g^{101}=g^{202}; \\
&g^{002} =-g^{112}, \quad
 g^{001} = -g^{212}.
\end{align}
It is then not difficult to check that these lead to \eqref{FB001}.

\subsection{The null condition in \cite{CLM18}}
Expanding \eqref{null3a}
in more details, we obtain
\begin{align}
N_{00} -(N_{01} +N_{10}) \cos \theta -(N_{02}+N_{20}) \sin \theta
+N_{11} \cos^2 \theta +N_{22} \sin^2 \theta + (N_{12}+N_{21}) \cos \theta \sin \theta =0.
\end{align}
Thus 
\begin{align}
& N_{00} +\frac 12 N_{11}+\frac 12 N_{22}=0; \\
& N_{01}+N_{10} =0, \quad N_{02}+N_{20} =0; \\
& N_{12}+N_{21} =0, \quad  N_{11} -N_{22}=0.
\end{align}
Therefore $N_{11}=N_{22}=-N_{00}$, and 
\begin{align}
N_{ij} a_i b_j = N_{00} (a_0b_0 -a_1b_1-a_2b_2)
+N_{01} (a_0b_1-a_1b_0) + N_{02} (a_0b_2-a_2b_0).
\end{align}
If $a=b$, then clearly $
N_{ij} a_i a_j = N_{00} (a_0^2 -a_1^2 -a_2^2).$
In particular,
\begin{align}
N_{ij} \partial_i v \partial_j v=N_{00} \cdot \Bigl( (\partial_t v)^2 -|\nabla v|^2 \Bigr).
\end{align}

\section{Proof of Theorem \ref{thm:main} }
In this section and Section 5, we carry out the proof of Theorem \ref{thm:main}. 
Write $v= \Gamma^{\alpha} u$, by Lemma \ref{Lem2.3} we have
\begin{align} 
\square v  &= \sum_{\alpha_1+ \alpha_2 \le \alpha}
g^{kij}_{\alpha_1,\alpha_2} \partial_k \Gamma^{\alpha_1} u 
\partial_{ij} \Gamma^{\alpha_2} u \label{3.0A}\\
&= 
g^{kij}\partial_k v
\partial_{ij} u  + 
g^{kij} \partial_k  u 
\partial_{ij} v 
+\sum_{\substack{\alpha_1<\alpha, \alpha_2<\alpha;\\ \alpha_1+ \alpha_2 \le \alpha}}
g^{kij}_{\alpha_1,\alpha_2} \partial_k \Gamma^{\alpha_1} u 
\partial_{ij} \Gamma^{\alpha_2} u. 
\end{align}

Choose $p(t,r) = q(r-t)$, where
\begin{align}
q(s) = \int_0^s \langle \tau\rangle^{-1}  \bigl(\log ( 2+\tau^2) \bigr)^{-2} d\tau, \quad s\in \mathbb R.
\end{align}
Clearly 
\begin{align}
-\partial_t p = \partial_r p = q^{\prime}(r-t)
= \langle r -t \rangle^{-1} \bigl( \log (2+(r-t)^2)  \bigr)^{-2}.
\end{align}
Multiplying both sides of \eqref{3.0A} by $e^p \partial_t v$, we obtain
\begin{align*}
   \text{LHS}&=\int e^p \pa_{tt}v\pa_{t}v -\int e^p \Delta v\pa_{t}v
   =\int e^p \pa_{tt} v\pa_{t}v +\int e^p\na v\cdot\na\pa_{t}v+\int e^p\na v\cdot\na p\pa_{t}v \\
   &=\frac12\frac{d}{dt}\int e^p(\pa v)^2 -\frac12\int e^{p}|\pa v|^2 p_{t}+\int e^p\na v\cdot\na p\pa_{t}v \\
   &=\frac12\frac{d}{dt}\| e^\frac{p}{2}\pa v\|_{L^2}^2 +\frac 12 
   \int e^p q^{\prime} \cdot \Bigl( |\partial_+ v|^2+ \frac {|\partial_{\theta} v|^2}{r^2}  \Bigr)
   =\frac12\frac{d}{dt}\| e^\frac{p}{2}\pa v\|_{L^2}^2 +\frac 12 
   \int e^p q^{\prime} |T v|^2.
\end{align*}

To simplify the notation in the subsequent nonlinear estimates, we introduce the following
terminology.

\noindent
\textbf{Notation}.  For a quantity $X(t)$, we shall write $X(t)= \mathrm{OK}$  if $X(t)$ can be written as
\begin{align}
X(t)= \frac d {dt} X_1 (t)+X_2(t) +X_3(t), 
\end{align}
where  (below $\alpha_0>0$ is some constant)
\begin{align}
|X_1(t)| \ll \| (\partial \Gamma^{\le m}  u)(t,\cdot) \|_{L_x^2(\mathbb R^2)}^2, 
\quad |X_2(t)| \ll    \sum_{|\alpha| \le m} \int e^p q^{\prime} |(T \Gamma^{\alpha} u)(t,x) |^2 dx,
\quad |X_3(t)| \lesssim \langle t \rangle^{-1-\alpha_0}.
\end{align}
In yet other words, the quantity $X$ will be controllable  if either it can be absorbed into
the energy, or can be controlled by the weighted $L^2$-norm of the good
unknowns  from the Alinhac weight,  or it is integrable in time. 
%It should
%be noted that here in the notation  ``$\ll$" the smallness of the pre-factors  
%is taken for $t$ sufficiently large. For example, the notation
%$|X_1(t)| \ll  \| (\partial \Gamma^{\le m} u ) (t,\cdot) \|_{L_x^2(\mathbb R^2)}^2$
%means for any $\epsilon>0$, we can find $T_0>0$ such that
%\begin{align}
%|X_1(t)| \le \epsilon \| (\partial \Gamma^{\le m}  u)(t,\cdot) \|_{L_x^2(\mathbb R^2)}^2,
%\qquad \forall\, t\ge T_0.
%\end{align}

We now proceed with the nonlinear estimates. We shall discuss several cases.

\subsection{The case $\alpha_1<\alpha$ and $\alpha_2<\alpha$}
Since $g^{kij}_{\alpha_1,\alpha_2}$ still
satisfies the null condition, by \eqref{2.10A} we have
\begin{align}
  &\sum_{\substack{\alpha_1<\alpha,  \alpha_2 <\alpha \\ \alpha_1+\alpha_2\le \alpha}}
g^{kij}_{\alpha_1,\alpha_2} \partial_k \Gamma^{\alpha_1} u 
\partial_{ij} \Gamma^{\alpha_2} u    \notag \\
= &\sum_{\substack{\alpha_1<\alpha,  \alpha_2 <\alpha \\ \alpha_1+\alpha_2\le \alpha}}
g^{kij}_{\alpha_1,\alpha_2}
(T_k \Gamma^{\alpha_1} u \partial_{ij} \Gamma^{\alpha_2} u
-\omega_k \partial_t \Gamma^{\alpha_1 } u T_i \partial_j \Gamma^{\alpha_2} u
+ \omega_k \omega_i \partial_t \Gamma^{\alpha_1} u T_j
\partial_t \Gamma^{\alpha_2} u ). 
\end{align}

\texttt{Estimate of $\| T_k \Gamma^{\alpha_1} u  \partial^2 \Gamma^{\alpha_2} u\|_2$}.
If $|\alpha_1|\le |\alpha_2|$,  then by Lemma \ref{lem2.6} we have
\begin{align}
 \| T_k \Gamma^{\alpha_1} u  \partial^2 \Gamma^{\alpha_2} u\|_2
 \lesssim
 \| \frac {T_k \Gamma^{\alpha_1} u } {\langle r-t\rangle}
 \|_{\infty} \cdot
 \| \langle r-t\rangle \partial^2 \Gamma^{\alpha_2}u \|_2 
 \lesssim  t^{-\frac 32}.
 \end{align}
 If $|\alpha_1|>|\alpha_2|$, then we have
 \begin{align}
 \| T_k \Gamma^{\alpha_1} u  \partial^2 \Gamma^{\alpha_2} u\|_2
 \lesssim
 \| \frac {T_k \Gamma^{\alpha_1} u } {\langle r-t\rangle}
 \|_{2} \cdot
 \| \langle r-t\rangle \partial^2 \Gamma^{\alpha_2}u \|_{\infty}
 \lesssim t^{-1} \cdot  t^{-\frac 12} \log t
 \lesssim t^{-\frac 32} \log t.
 \end{align}

\texttt{Estimate of $\| \partial \Gamma^{\alpha_1} u T \partial \Gamma^{\alpha_2} u\|_2$}.
If $|\alpha_1| \le |\alpha_2|$ we have
\begin{align}
\| \partial \Gamma^{\alpha_1} u T \partial \Gamma^{\alpha_2} u\|_2
\lesssim \| \partial \Gamma^{\alpha_1} u \|_{\infty}
\cdot \| T \partial \Gamma^{\alpha_2} u \|_2 \lesssim t^{-\frac 32}.
\end{align}
If $|\alpha_1| > |\alpha_2|$ we have
\begin{align}
\| \partial \Gamma^{\alpha_1} u T \partial \Gamma^{\alpha_2} u\|_2
\lesssim \| \partial \Gamma^{\alpha_1} u \|_{2}
\cdot \| T \partial \Gamma^{\alpha_2} u \|_{\infty} \lesssim  t^{-\frac 32} \log t.
\end{align}

Collecting the estimates, we have proved
\begin{align} \label{4.13U}
  &\| \sum_{\substack{\alpha_1<\alpha,  \alpha_2 <\alpha \\ \alpha_1+\alpha_2\le \alpha}}
g^{kij}_{\alpha_1,\alpha_2} \partial_k \Gamma^{\alpha_1} u 
\partial_{ij} \Gamma^{\alpha_2} u    \|_2 \lesssim  t^{-\frac 32}\log t.
\end{align}

\subsection{The case $\al_{2}=\al$.}
Noting that $g^{kij}_{0,\alpha}=g^{kij}$, we have
\begin{align}
 \int g^{kij}\pa_{k}u\pa_{ij}v\pa_{t}v e^{p}
 &= \OK \underbrace{- \int g^{kij}\pa_{jk}u\pa_{i}v\pa_{t}v e^{p}}_{I_1}\underbrace{-\int g^{kij}\pa_{k}u\pa_{i}v\pa_{t}v \pa_{j}(e^{p})}_{I_2}-\int g^{kij}\pa_{k}u\pa_{i}v\pa_{tj}v e^{p}.
\end{align}
Here in the above, the term ``OK" is zero if $\partial_j =\partial_1$ or $\partial_2$. 
This term is nonzero when $\partial_j = \partial_t$, i.e. we should absorb it into the energy when
integrating by parts in the time variable. 

Further integration by parts gives
\begin{align}
  -\int g^{kij}\pa_{k}u\pa_{i}v\pa_{tj}v e^{p}
  &=\OK+ \underbrace{\int g^{kij}\pa_{tk}u\pa_{i}v\pa_{j}v e^{p}}_{I_3}+\underbrace{\int g^{kij}\pa_{k}u\pa_{i}v\pa_{j}v \pa_{t}(e^{p})}_{I_4}+\int g^{kij}\pa_{k}u\pa_{it}v\pa_{j}
  v e^{p}.
  \end{align}
\begin{align}
  \int g^{kij}\pa_{k}u\pa_{it}v\pa_{j}v e^{p}&=\OK \underbrace{-\int g^{kij}\pa_{ik}u\pa_{t}v\pa_{j}v e^{p}}_{I_5}\underbrace{-\int g^{kij}\pa_{k}u\pa_{t}v\pa_{j}v \pa_{i}(e^{p})}_{I_6}-\int g^{kij}\pa_{k}u\pa_{t}v\pa_{ij}v e^{p}.
\end{align}
It follows that
$$2\int g^{kij}\pa_{k}u\pa_{ij}v\pa_{t}v e^{p}=(I_1+I_3+I_5)+(I_2+I_4+I_6) +\OK.$$

Observe that if $\varphi=\pa_{k}u$ or $\varphi=e^{p}$, then 
\begin{align}\label{varphi1}
   &-\pa_{j}\varphi\pa_{i}v\pa_{t}v+\pa_{t}\varphi\pa_{i}v\pa_{j}v-\pa_{i}\varphi\pa_{t}v\pa_{j}v
   \notag\\
  =&-T_{j}\varphi\pa_{i}v\pa_{t}v+\omega_{j}\pa_{t}\varphi\pa_{i}v \pa_{t} v+\pa_{t}\varphi\pa_{i}v\pa_{j}v-T_{i}\varphi\pa_{t}v\pa_{j}v+\omega_{i}\pa_{t}\varphi\pa_{t}v T_{j}v-\omega_{i}\omega_{j}\pa_{t}\varphi(\pa_{t}v)^2 \notag \\
  =&-T_{j}\varphi\pa_{i}v\pa_{t}v+\pa_{t}\varphi\pa_{i}v T_{j} v-T_{i}\varphi\pa_{t}v\pa_{j}v+\omega_{i}\pa_{t}\varphi\pa_{t}v T_{j}v-\omega_{i}\omega_{j}\pa_{t}\varphi(\pa_{t}v)^2 \notag\\
  =&-T_{j}\varphi\pa_{i}v\pa_{t}v+\pa_{t}\varphi T_{i}v T_{j} v-T_{i}\varphi\pa_{t}v\pa_{j}v-\omega_{i}\omega_{j}\pa_{t}\varphi(\pa_{t}v)^2.
 \end{align}
By \eqref{varphi1} and rewriting $\partial_t \varphi = \partial_k \partial_t u
=T_k \partial_t u -\omega_k \partial_{tt} u$, we have
\begin{align}
I_1+I_3+I_5=&\;\int g^{kij}(
-T_{j}\partial_k u \pa_{i}v\pa_{t}v+\pa_{t}\partial_k u T_{i}v T_{j} v-T_{i}\partial_k u\pa_{t}v\pa_{j}v-\omega_{i}\omega_{j}T_k \partial_t u(\pa_{t}v)^2) e^pdx.
\end{align}
By Lemma \ref{lem2.6}, we have $\|T \partial u\|_{\infty} \lesssim t^{-\frac 32} \log t$
and $\| \langle r-t\rangle \partial^2 u \|_{\infty} \lesssim t^{-\frac 12} \log t$. Clearly then
\begin{align}
\int_{\text{$r<\frac t2$ or $r>2t$} } |\partial^2 u| |T v|^2 dx\lesssim t^{-\frac 32} \log t,
\quad \int_{r \sim t } |\partial^2 u | |Tv|^2 dx 
\ll \int e^p q^{\prime} |Tv|^2 dx. 
\end{align}
It follows that 
\begin{align}
I_1+I_3+I_5 =\OK.
\end{align}
Plugging $\varphi=e^{p}$ in \eqref{varphi1} and noting that $T_j (e^p)=0$, we  have
\begin{align*}
    I_2+I_4+I_6    
    &=\int g^{kij}\pa_k u\Big(-T_j(e^p)\pa_iv \pa_tv- T_i(e^p)\pa_tv\pa_jv-\omega_i\omega_j(\pa_tv)^2\pa_t(e^p)+\pa_t(e^p)T_iv T_jv \Big) \notag \\
  &=   \int g^{kij}\left(  -T_k u\cdot \omega_i\omega_j (\pa_tv)^2 \pa_t(e^p) +\partial_k u\pa_t(e^p)T_i v T_j v \right).
\end{align*}
By Lemma \ref{lem2.6} we have $||Tu| \partial_t (e^p)|\lesssim t^{-\frac 32}$. 
Clearly
\begin{align}
\|\partial u \partial_t (e^p)\|_{L_x^{\infty}(r<\frac t2,\, \text{or } r>2t)} \lesssim t^{-\frac 32},
\quad \int_{r\sim t} |\partial u \partial_t (e^p)| |Tv|^2 dx \ll \int e^p q^{\prime} |Tv|^2 dx. 
\end{align}
Thus
\begin{align*}
  I_{2}+I_{4}+I_{6} = \OK. 
\end{align*}
This concludes the case $\alpha_2=\alpha$. In the next section we deal with the main
piece $\alpha_1=\alpha$. 

\section{Estimate of the main piece $\al_{1}=\al$, $\al_{2}=0$}
In this section we estimate the main piece $\alpha_1=\alpha$. By \eqref{2.10A}, we have 
\begin{align*}
 \int g^{kij}\pa_{k}v\pa_{ij}u\pa_{t}v e^{p}
  =&\int g^{kij} (T_{k}v\pa_{ij}u-\omega_{k}\pa_{t}vT_{i}\pa_{j}u+\omega_{k}\omega_{i}\pa_{t}v T_{j}\pa_{t}u)\pa_{t}v e^{p}\\
  =&\int g^{kij} (T_{k}vT_{i}\pa_{j}u-\omega_{i}T_{k}vT_{j}\pa_{t}u+\omega_{i}\omega_{j}T_{k}v\pa_{tt}u-\omega_{k}\pa_{t}vT_{i}\pa_{j}u+\omega_{k}\omega_{i}\pa_{t}v T_{j}\pa_{t}u)\pa_{t}v e^{p}.
\end{align*}
By Lemma \ref{lem2.6}, all terms containing $T\partial u$ decay as $O(t^{-\frac 32} \log t)$.
Thus 
\begin{equation}\label{eq:al1-00}
  \int g^{kij} \pa_{k}v\pa_{ij}u\pa_{t}v e^{p}= \OK+ \int g^{kij} \omega_{i}\omega_{j}T_{k}v\pa_{tt}u\pa_{t}v e^{p}. 
\end{equation}
Recall $T_{0}=0$, $T_{1}=\omega_{1}\pa_{+}-\frac{\omega_{2}}{r}\pa_{\theta}$,  $T_{2}=\omega_{2}\pa_{+}+\frac{\omega_{1}}{r}\pa_{\theta}$. By \eqref{null2},   we have
\begin{align*}
 g^{kij} \omega_{i}\omega_{j}T_{k}v
 &= g^{1ij} \omega_{i}\omega_{j}(\omega_{1}\pa_{+}v-\frac{\omega_{2}}{r}\pa_{\theta}v)
  +g^{2ij} \omega_{i}\omega_{j}(\omega_{2}\pa_{+}v+\frac{\omega_{1}}{r}\pa_{\theta}v)\\
 &=\left(g^{1ij} \omega_{1}\omega_{i}\omega_{j}+g^{2ij} \omega_{2}\omega_{i}\omega_{j}\right)\pa_{+}v\coloneqq h(\theta)\pa_{+}v.
\end{align*}
Choose nonnegative $\phi \in C_c^{\infty}(\mathbb R)$ such that $\phi(z)=1$ for 
$\frac 23 \le z \le \frac 32$ and $\phi(z)=0$ for $z\le \frac 13$ or $z\ge 2$. Then
\begin{equation}\label{eq:al1-01}
\begin{split}
  \int g^{kij} \omega_{i}\omega_{j}T_{k}v\pa_{tt}u\pa_{t}v e^{p}
  &= \int h(\theta)\pa_{+}v\pa_{tt}u\pa_{t}v e^{p}\cdot\left(1-\phi\left(\frac rt\right)\right)+\int h(\theta)\pa_{+}v\pa_{tt}u\pa_{t}v e^{p}\phi\left(\frac rt\right).
  \end{split}
\end{equation}
By Lemma \ref{lem2.6}, we have
\begin{equation}\label{eq:al1-02}
  \int h(\theta)\pa_{+}v\pa_{tt}u\pa_{t}v e^{p}\cdot\left(1-\phi\right)\lesssim t^{-1}\int |\pa v|^2|\Lg r-t\Rg\pa_{tt}u| \lesssim t^{-\frac 32} \log t = \OK.
\end{equation}
By using the identity $\pa_{t}=\frac{\pa_{+}+\pa_{-}}{2}$ and 
the fact that $\| \langle r -t \rangle \partial^2 u \|_{\infty} \lesssim t^{-\frac 12} \log t$, we get
\begin{align}\label{eq:al1-03}
  2\int h(\theta)\pa_{+}v\pa_{tt}u\pa_{t}v e^{p}\phi &=\int h(\theta)\pa_{+}v\pa_{tt}u\pa_{+}v e^{p}\phi+\int h(\theta)\pa_{+}v\pa_{tt}u\pa_{-}v e^{p}\phi \notag \\
  &=\OK+ \int h(\theta)\pa_{+}v\pa_{tt}u\pa_{-}v e^{p}\phi.
\end{align}
Integrating by parts, we have
\begin{align}
  \int h(\theta)\pa_{+}v\pa_{tt}u\pa_{-}v e^{p}\phi  \cdot  rdr d\theta
  &=\;\frac{d}{dt}\int h(\theta)v\pa_{tt}u\pa_{-}v e^{p}\phi dx
   -\int h(\theta)v\pa_{-}v \pa_{+}\left(\pa_{tt}ue^{p}\phi\right) dx
  \notag \\
  & \quad -\int h(\theta)v\pa_{tt}u\pa_{+}\pa_{-}v e^{p}\phi dx
  -\int h(\theta) v \partial_{tt} u \partial_- v e^p \phi \frac 1r dx.
 \end{align}
 In the above computation, one should note that when integrating by parts in $r$ we should
 take into consideration the metric $rdr$. The fourth term exactly corresponds to the derivative
 of the metric factor.
The first  and fourth terms are clearly acceptable by using Hardy and the decay of $\langle r-t\rangle \partial_{tt}u$. 
For the second term we have
\begin{align}
  \left| \langle r -t\rangle \pa_{+}\left(\pa_{tt}ue^{p}\phi\right)\right|
  &\lesssim \left|\langle r -t\rangle \pa_{+}\pa_{tt}u\phi\right|+\left|\langle r-t\rangle \pa_{tt}u\pa_{+}\phi\right|  \notag \\
  &\lesssim  t^{-1}\|2\langle r -t \rangle L_0\pa_{tt}u- \langle r-t\rangle ( t-r)\pa_{-}\pa_{tt}u
  \|_{L_x^{\infty}(|x|>\frac t{10})} + t^{-\frac 32}\lesssim t^{-\frac32}. \label{5.6A}
\end{align}
Here in the derivation of \eqref{5.6A}, we used Lemma \ref{lem2.6} and the inequalities
\begin{align}
&|\langle r -t\rangle L_0  \partial_{tt} u|
\lesssim | \langle r -t \rangle \partial_{tt} \Gamma^{\le 1} u |
\lesssim t^{-\frac 12}, \qquad \text{for }r \ge t/10.
\end{align}
For the third term we use the identity $\pa_{+}\pa_{-}v=\Box v+\frac{\pa_{r}v}{r}+\frac{\pa_{\theta\theta}v}{r^2}$ and compute it as
\begin{equation}\label{eq:al1-2}
 \begin{split}
  &\int h(\theta)v\pa_{tt}u\pa_{+}\pa_{-}v e^{p}\phi\\
  =&\int h(\theta)v\pa_{tt}u\left(\frac{\pa_{r}v}{r} +\frac{\pa_{\theta\theta}v}{r^2}\right)e^{p}\phi
  +\sum_{\be_{1}+\be_{2}\leq \al}\int h(\theta)v\pa_{tt}u\cdot g_{\be_{1},\be_{2}}^{kij}\pa_{k}\Ga^{\be_{1}}u\pa_{ij}\Ga^{\be_{2}}ue^{p}\phi.
 \end{split}
\end{equation}
Integrating by parts (for the term $\partial_{\theta\theta}v$), we have
\begin{align*}
  &\int h(\theta)v\pa_{tt}u\left(\frac{\pa_{r}v}{r}+\frac{\pa_{\theta\theta}v}{r^2}\right)e^{p}\phi\\
  =& \int h(\theta) \frac {v}{\langle r -t\rangle}
  \langle r -t\rangle \partial_{tt} u \partial_r v \cdot \frac 1r e^p \phi 
  -\int h(\theta)\pa_{tt}u\left(\frac{\pa_{\theta}v}{r}\right)^2e^{p}\phi-\int \pa_{\theta}(h(\theta)\pa_{tt}u)v\frac{\pa_{\theta}v}{r^2}e^{p}\phi\\
  =& \OK.
  \end{align*}
By \eqref{4.13U}, we have
\begin{equation*}
  \sum_{\be_{1}< \al,\be_{2}< \al,\atop \be_{1}+\be_{2}\leq \al}\int h(\theta)v\pa_{tt}u \cdot g_{\be_{1},\be_{2}}^{kij}\pa_{k}\Ga^{\be_{1}}u\pa_{ij}\Ga^{\be_{2}}u e^{p}\phi\lesssim   t^{-2}=\OK.
\end{equation*}
For the term  $\be_{1}=\al$, $\be_{2}=0$ in \eqref{eq:al1-2}, it follows from \eqref{a2.12a} that
\begin{align*}
  \int g ^{kij}h(\theta)v\pa_{tt}u \pa_{k}v\pa_{ij}u e^{p}\phi
  \lesssim &\int|v\pa_{tt}u||T v\pa^2u|e^{p}\phi+\int|v\pa_{tt}u||\pa v||T\pa u| e^{p}\phi\\
  \lesssim &\int|T v|^2|\pa^2u|e^{p}\phi+t^{-\frac32}\left\|\Lg r-t\Rg^{-1}v\right\|_{L_{x}^2(\R^2)}^2+t^{-2}\\
=& \OK.
\end{align*}

For the term  $\be_{1}=0$, $\be_{2}=\al$ in \eqref{eq:al1-2}, we apply \eqref{2.10A} to obtain
\begin{align*}
  \int g ^{kij}h(\theta)v\pa_{tt}u\pa_{k}u \pa_{ij}v e^{p}\phi
  = \int g ^{kij}h(\theta)v\pa_{tt}u\cdot(T_{k}u\pa_{ij}v-\omega_{k}\pa_{t}uT_{i}\pa_{j}v+\omega_{k}\omega_{i}\pa_{t}u T_{j}\pa_{t}v) e^{p}\phi.
\end{align*}
We rewrite it as
\begin{align*}
  \int g ^{kij}h(\theta)v\pa_{tt}uT_{k}u\pa_{ij}v e^{p}\phi
  =&\int g ^{kij}\pa_{i}(h(\theta)v\pa_{tt}uT_{k}u\pa_{j}v e^{p}\phi)
  -\int g ^{kij}\pa_{i}(h(\theta) T_{k}u e^{p}\phi )v \pa_{tt}u\pa_{j}v\\
  &-\int g ^{kij}h(\theta)v\pa_{i}\pa_{tt}uT_{k}u\pa_{j}v e^{p}\phi-\int g ^{kij}h(\theta)\pa_{i}v\pa_{tt}uT_{k}u\pa_{j}v e^{p}\phi.
\end{align*}
The term $\int g ^{kij}\pa_{i}(h(\theta)v\pa_{tt}uT_{k}u\pa_{j}v e^{p}\phi)$ is 
zero for $i\ne 0$. For $i=0$ it is clearly acceptable since it can be absorbed into the time
derivative of the energy due to its smallness. By Lemma \ref{lem2.3a} and \ref{lem2.6}, we have
\begin{align*}
  |\pa_{i}(h(\theta) T_{k}u e^{p}\phi )|&\lesssim |\pa_{i}h(\theta) T_{k}u e^{p}\phi |+|h(\theta)\pa_{i} T_{k}u e^{p}\phi |+|h(\theta) T_{k}u \pa_{i}e^{p}\phi |+|h(\theta) T_{k}u e^{p}\pa_{i}\phi |\\
  &\lesssim  t^{-\frac32}+|h(\theta)\pa_{i} \omega_{k}\pa_{t}u e^{p}\phi |+|h(\theta)T_{k}\pa_{i}u e^{p}\phi |+\left|h(\theta) \frac{T_{k}u}{\Lg r-t\Rg} \phi \right|
  \lesssim  t^{-\frac32}.
 \end{align*}
The term containing $v\partial_i \partial_{tt} u$ can be handled by \eqref{2.30D}. Thus
\begin{align*}
  &\int g ^{kij}h(\theta)v\pa_{tt}uT_{k}u\pa_{ij}v e^{p}\phi
  = \OK.
\end{align*}
Similarly, we have
\begin{align*}
  &\int g ^{kij}\omega_{k}h(\theta)v\pa_{tt}u\pa_{t}uT_{i}\pa_{j}v e^{p}\phi
  =\int g ^{kij}\omega_{k}h(\theta)v\pa_{tt}u\pa_{t}u\left(\pa_{j}T_{i}v-\pa_{j}\omega_{i}\pa_{t}v \right) e^{p}\phi =\OK,\\
  &\int g ^{kij}\omega_{k}\omega_{i}h(\theta)v\pa_{tt}u\pa_{t}uT_{j}\pa_{t}v e^{p}\phi=\int g ^{kij}\omega_{k}\omega_{i}h(\theta)v\pa_{tt}u\pa_{t}u\pa_{t}T_{j}v e^{p}\phi = \OK.
\end{align*}
This concludes the estimate of the main piece.

\bibliographystyle{abbrv}

%\appendix

\end{document}